\documentclass[12pt,a4paper]{amsart}
\usepackage{amscd}
\theoremstyle{plain} %% This is the default
\newtheorem{thm}{Theorem}[section]
\newtheorem{cor}[thm]{Corollary}
\newtheorem{lem}[thm]{Lemma}
\newtheorem{prop}[thm]{Proposition}

\theoremstyle{definition}

\newtheorem{rem}[thm]{Remark}

\theoremstyle{remark}

\numberwithin{equation}{section}

\newsavebox{\SmallMathBox}

\def\dd{\partial}

\def\Di{D\kern -.65em /}
\def\Dii{D\kern -.45em /}
\def\di{{\dd}\kern -.55em /}
\def\dii{{\dd}\kern -.40em /}

\def\noi{\noindent}

\def\tand{\mbox{\ \rm  and }}

\def\ZZ{{\bf Z}}

\def\Cc{{\mathcal C}}
\def\Dd{{\mathcal D}}
\def\Ee{{\mathcal E}}
\def\Ff{{\mathcal F}}
\def\Gg{{\mathcal G}}

\def\Kk{{\mathcal K}}

\def\Mm{{\mathcal M}}
\def\Nn{{\mathcal N}}

\def\Rr{{\mathcal R}}
\def\Ss{{\mathcal S}}
\def\Tt{{\mathcal T}}

\def\={\cong}
\def\>{\supset}
\def\<{\subset}

\def\12{\frac{1}{2}}
\def\2{\Dd}
\def\3{\Nn}
\def\4{\Rr}
\def\6{\cup}
\def\8{\otimes}
\def\0{^{\circ}}

\def\a{\alpha}

\def\e{\varepsilon}
\def\g{\gamma}
\def\G{\Gamma}

\def\k{\kappa}

\def\la{\lambda}

\def\p{\pi}

\def\s{\sigma}
\def\Si{\Sigma}
\def\z{\zeta}
\def\Z{\ZZ}

\def\Si{S\kern -.65em /}

\def\tensor{\otimes}

\def\over{/}

\begin{document}

\title[Adiabatic Decomposition of the $\z$-determinant I.]
{Adiabatic Decomposition of the $\z$-determinant of the Dirac Laplacian I.\\
          The Case of an Invertible Tangential Operator}

\author[Jinsung Park]{Jinsung Park*}
\address{Department of Mathematics\\IUPUI (Indiana/Purdue)\\
Indianapolis IN 46202--3216, U.S.A.} \email{jinspark@indiana.edu}

\author{Krzysztof P. Wojciechowski}
\address{Department of Mathematics\\IUPUI (Indiana/Purdue)\\
Indianapolis IN 46202--3216, U.S.A.}
\email{kwojciechowski@math.iupui.edu}

\address{Department of Mathematics, Inha University, Inchon,
402-751, Korea.} \email{ywlee@math.inha.ac.kr}

\thanks{*Partially supported by Korea Science and Engineering
Foundation}

%\date{{\em \today. File name:} CPDE7.tex}

\maketitle

\centerline{with an Appendix by Yoonweon Lee}

\vskip 5mm

\begin{abstract}
We discuss the decomposition of the $\z$-determinant of the square
of the Dirac operator into the contributions coming from the
different parts of the manifold. The result was announced in
\cite{JPKW1}\,. The proof sketched in \cite{JPKW1} was based on
results of Br${\ddot {\rm u}}$ning and Lesch (see \cite{BL98}). In
the meantime we have found another proof, more direct and
elementary, and closer to the spirit of the original papers which
initiated the study of the adiabatic decomposition of the spectral
invariants (see \cite{DoWo91} and \cite{Si88}). We discuss this
proof in detail. We study the general case (non-invertible
tangential operator) in forthcoming work (see \cite{JKPW2} and
\cite{JKPW3}). In the Appendix we present the computation of the
cylinder contribution to the $\z$-function of the Dirac Laplacian
on a manifold with boundary, which we need in the main body of the
paper. This computation is also used to show the vanishing result
for the $\z$-function on a manifold with boundary.
\end{abstract}

%\maketitle

\section*{Results}

Let $\Dd : C^{\infty}(M;S) \to C^{\infty}(M;S)$ be a compatible
Dirac operator acting on sections of a bundle of Clifford modules
$S$ over a closed manifold $M$. Assume that we have a
decomposition of $M$ as $M_1 \cup M_2$\,, where $M_1$ and $M_2$
are compact manifolds with boundary such that
\begin{equation}\label{e:dec1}
M = M_1 \cup M_2 \ \ , \ \ M_1 \cap M_2 = Y = {\partial}M_1 =
{\partial}M_2 \, .
\end{equation}
The $\z$-determinant of the operator $\Dd$ is given by the formula
\begin{equation}\label{e:zd}
det_{\z}\Dd = e^{{\frac{i\p}{2}}(\z_{\Dd^2}(0) -
\eta_{\Dd}(0))}{\cdot} e^{-{\frac1{2}}\z'_{\Dd^2}(0)} \,,
\end{equation}
(see \cite{Si85}, see also the Introduction of \cite{SSKPW299}).
In this paper we study the decomposition of $det_{\z}\Dd$ on $M$
into contributions coming from $M_1$ and $M_2$\,. This issue was
already solved for the phase of the determinant
$${\frac{i\p}{2}}(\z_{\Dd^2}(0) - \eta_{\Dd}(0)) \,,$$
and there remains only the modulus - the {\it square root} of the
$\z$-determinant of the Dirac Laplacian $\Dd^2$ - to study. We
present here an {\it ``adiabatic"} solution of the problem in the
case of an {\it ``invertible tangential operator"}. The general
case will be presented in \cite{JKPW3} (see also \cite{JKPW2}).
However, the discussion in this paper is an important part of the
study of the general case.

\bigskip

We start with a brief discussion of the splitting of the phase of
the $\z$-determinant. The invariant $\z_{\Dd^2}(0)$ poses no
problems. The value of the function $\z_{\Dd^2}(s)$ at $s=0$ is a
local invariant in the sense that it is given by a formula
$$\z_{\Dd^2}(0) = \int_M a(x)dx \,,$$
where $a(x)$ is a density determined at the point $x \in M$ by the
coefficients of the operator $\Dd$ at this point (see for instance
\cite{Gi95}). This is the reason why the index of an elliptic
differential operator, which can  be viewed as the difference of
the values of two different $\z$-functions determined by the
operator $\Dd$, has a nice  decomposition corresponding to the
decomposition of the manifold.

The other contribution to the phase of $det_{\z}\Dd$ is the
eta-invariant $\eta_{\Dd}(0)$ and this is not a local invariant
(see \cite{AtPaSi75a}), hence at first sight it is difficult to
expect a nice and clear splitting formula. It is therefore rather
surprising that such a formula for $\eta_{\Dd}(0)$ actually
exists.

In the following we concentrate on the odd-dimensional case
$$n = dim \ M = 2k+1 \,. $$
We further assume that $M$ and the operator $\Dd$ have product
structures in a neighborhood of the boundary $Y$. More precisely,
we assume that there is a bicollar neighborhood $N = [-1,1] \times
Y$  of $Y$ in $M$ such that  the Riemannian structure on $M$ and
the Hermitian structure on $S$ are products when restricted to
$N$. This implies that $\Dd$ has the following form when
restricted to the submanifold $N$
\begin{equation}\label{e:pr}
\Dd = G(\partial_u + B) \, .
\end{equation}
Here $u$ denotes the normal variable, $G : S|_Y \to S|_Y$ is a
bundle automorphism, and $B$ is a corresponding Dirac operator on
$Y$. Moreover, $G$ and $B$ do not depend on $u$ and they satisfy
\begin{equation}\label{e:bd}
G^* = - G \ \ , \ \  G^2 = -Id \ , \ \ B = B^* \ \ \tand \ \ GB =
- BG \ \, .
\end{equation}
The operator $B$ has a discrete spectrum with infinitely many
positive and infinitely many negative eigenvalues. In this work we
consider only the case of an invertible tangential operator, i.e.
we assume that $ker \ B = \{0\}$\,. The general case is more
difficult to handle and we refer to \cite{JKPW2} and \cite{JKPW3}
for the discussion of the {\it noninvertible case}. However, the
present work plays an important part in the analysis of the
general case.

\bigskip

Let $\Pi_>$ denote the spectral projection onto the subspace
spanned by the eigensections of $B$ corresponding to the positive
eigenvalues. Then $\Pi_>$ is an elliptic boundary condition for
$\Dd_2 = \Dd|_{M_2}$ (see \cite{AtPaSi75}; see \cite{BoWo93} for
an exposition of the theory of elliptic boundary problems for
Dirac operators). In fact, any orthogonal projection satisfying
\begin{equation}\label{e:gr*}
-GPG = Id - P \ \ \tand \ \ P - \Pi_> \ \ \text{{\it is  a
smoothing operator}},
\end{equation}
is a self-adjoint elliptic boundary condition for the operator
$\Dd_2$. This means that the associated operator
$$(\Dd_2)_P : dom \ (\Dd_2)_P \to
L^2(M_2;S|_{M_2})$$ with $dom \ (\Dd_2)_P = \{s \in
H^1(M_2;S|_{M_2}) \mid P(s|_Y) = 0\}$  is a self-adjoint Fredholm
operator with $ker((\Dd_2)_P)\subset C^{\infty}(M_2;S|_{M_2})$ and
a discrete spectrum (see \cite{KPW99}).

The existence of the meromorphic extensions of the functions
$\eta_{(\Dd_2)_P}(s)$, $\z_{(\Dd_2)_P^2}(s)$ to the whole complex
plane and their nice behavior in a neighborhood of $s=0$ was
established in \cite{KPW99}. We denote by
${\Gg}r^*_{\infty}(\Dd_2)$ the space of P satisfying
(\ref{e:gr*}).

Let us observe that $Id - P \in {\Gg}r^*_{\infty}(\Dd_1)$\,, if
$P$ is an element of ${\Gg}r^*_{\infty}(\Dd_2)$. We denote by
$\eta_{G(\partial_u + B)}(P_1, P_2)(s)$ the $\eta$-function of the
operator $G(\partial_u + B)$ on $[0,1] \times Y$ subject to the
boundary condition $P_2$ at $u=0$ and $Id-P_1$ at $u=1$\,. We have
the following pasting formula proved in \cite{KPW99}
\begin{equation}\label{e:split12}
\eta_{\Dd}(0) = \eta_{(\Dd_{1})_{Id-P_1}}(0) +
\eta_{(\Dd_{2})_{P_2}}(0) + \eta_{G(\partial_u + B)}(P_1, P_2)(0)
\ \ mod \ \Z \,.
\end{equation}

\noi A similar formula for finite-dimensional perturbations of
$\Pi_>$ has been discussed by several authors (see \cite{KPW94,
KPW95, KPW99} and references therein).

\bigskip
\noi The proof of (\ref{e:split12}) offered by the second author
goes as follows.

\noi First, we replace the bicollar $N$ by $N_R = [-R,R] \times
Y$. Now $\eta_{\Dd}(0)$, which can be expressed using an
appropriate heat-kernel formula, splits into contributions coming
from each side, plus the cylinder contribution (vanishing in the
case of $\Dd$) and error terms. The error terms disappear as $R
\to \infty$\,.

\noi Second, though $\eta_{\Dd}(0)$ is not local, its variation
(for instance with respect to the parameter $R$) is local and
therefore the value of the contributions does not vary with $R$.
This is enough to make explicit calculations of the formula
(\ref{e:split12}).

\bigskip

In this work we apply the strategy employed above to study
$$det_{\z}\Dd^2 =
e^{-{\frac d{ds}}\z_{\Dd^2}(s)|_{s=0}} \,.$$ However, we have to
take into account two additional difficulties, which arise in the
case of the $\z$-determinant of $\Dd^2$\,.

First of all, the invariant $-{\frac d{ds}}\z_{\Dd^2}(s)|_{s=0}$
is much more subtle than the $\eta$-invariant. Even the variation
of $-{\frac d{ds}}\z_{\Dd^2}(s)|_{s=0}$ is not given by a local
formula.

Second, the cylinder contribution is not trivial in this case.

We handled those difficulties in \cite{JPKW1} using the technique
developed in \cite{BL98}. Here we choose a different path. The
invariant ${\frac d{ds}}\z_{\Dd^2}(s)|_{s=0}$ is given by the
formula
\begin{equation}\label{e:dt1}
{\frac d{ds}}\z_{\Dd^2}(s)|_{s=0} = \int_0^{\infty} {\frac 1{t}}Tr
\ e^{-t\Dd^2} dt \, .
\end{equation}
Let us explain how to interpret formula (\ref{e:dt1}). The trace
$Tr \ e^{-t\Dd^2}$ has an asymptotic expansion of the form
$$Tr \ e^{-t\Dd^2} =
t^{-{\frac n{2}}}\sum_{k=0}^N a_kt^k + O(t^{N + {\frac {1-n}{2}}})
\,,$$ where $a_k = \int_M\a_k(x)dx$\,, and the density $\a_k(x)$
at the point $x \in M$ is determined by coefficients of the
operator $\Dd^2$ (see \cite{Gi95}). This shows that
$$\z_{\Dd^2}(s) =
{\frac1{\G(s)}}\int_0^{\infty}t^{s-1} Tr \ e^{-t\Dd^2}dt$$

\noi is a holomorphic function of $s$\,, for $Re(s) > {\frac
n{2}}$\,, and that it has a meromorphic extension to the whole
complex plane with (possible) simple poles at $s_k = {\frac n{2}}
- k$\,. The $\G$-function has the following form in a neighborhood
of $0$
$$\G(s) =
{\frac 1{s}} + \g + s\,h(s) \,,$$ where $\g$ is Euler's constant
and $h(s)$ is a holomorphic function in a neighborhood of $0$\,.
This allows us to compute $\z_{\Dd^2}(0)$
\begin{align*}
\z_{\Dd^2}(0) &= \lim_{s \to 0} {\frac
1{\G(s)}}\int_0^{\infty}t^{s-1} Tr \ e^{-t\Dd^2}dt = \lim_{s \to
0}
s\int_0^1t^{s-1} Tr \ e^{-t\Dd^2}dt\\
&= \lim_{s \to 0} s\int_0^1t^{s-1}t^{-{\frac n{2}}}(\sum_{k=0}^N
t^ka_k)dt = \lim_{s \to 0}s{\cdot}\sum_{k=0}^N {\frac{2a_k}{2s +
2k - n}} = a_{n\over{2}} \,,
\end{align*}
where $N$ denotes any sufficiently large natural number and we
keep in mind that

$$a_{n\over{2}} = 0 \ \ \text{for $n$
odd}.$$

\bigskip
\noi In particular, $\z_{\Dd^2}(0) = 0$ for $n$ odd. Though $s=0$
is a regular point, the $\z$-function may have poles on the right
side of $0$\,, and the function

$$\k_{\Dd^2}(s) =
\int_0^{\infty}t^{s-1} Tr \ e^{-t\Dd^2}dt$$

\noi has even more poles. In particular, following the
computations presented above, we have
$$Res_{s=0}\k_{\Dd^2}(s) = a_{n\over{2}} \,.$$

\bigskip
\noi Now, the derivative of the $\z$-function at $s=0$ is obtained
as follows

$$\z_{\Dd^2}'(0) = {\frac
d{ds}}({\frac{\k_{\Dd^2}(s)}{\G(s)}})|_{s=0} = {\frac
d{ds}}({\frac{a_{n\over{2}} + s(\k_{\Dd^2}(s) -{\frac
{a_{n\over{2}}}{s}})}{1 + s\g + s^2h(s)}})|_{s=0} =$$

\bigskip

$${\frac{(\k_{\Dd^2}(s) -{\frac{a_{n\over{2}}}{s}})(1 + s\g
+ s^2h(s)) - (a_{n\over{2}} + s(\k_{\Dd^2}(s) -
{\frac{a_{n\over{2}}}{s}}))(\g + 2sh(s))}{(1 + s\g +
s^2h(s))^2}}|_{s=0} =$$

\bigskip

$$(\k_{\Dd^2}(s) -{\frac{a_{n\over{2}}}{s}})|_{s=0} -
{\g}a_{n\over{2}} = (\k_{\Dd^2}(s)
-{\frac{a_{n\over{2}}}{s}})|_{s=0} - {\g}a_{n\over{2}} \ \ .$$

\bigskip
\noi This discussion provides a justification for the a priori
{\it ``formal"} formula (\ref{e:dt1}).

\bigskip

\begin{rem}\label{r:z0}

(a) For simplicity we presented here the $\z$-function in the case
$ker \ \Dd = \{0\}$ . In general we define $\z$-function as

$$\z_{\Dd^2}(s) =
{\frac1{\G(s)}}\int_0^{\infty}t^{s-1} (Tr \ e^{-t\Dd^2} - dim \
ker \ \Dd)dt \ \ ,$$

and

$$\z_{\Dd^2}(0) = a_{n\over{2}} - dim \ ker \ \Dd \ \
.$$

\bigskip

\noi (b) The corresponding result for the boundary value problems
is proved in the Appendix (see also \cite{YLKW01}). It is shown
that
$$\z_{\Dd_{i,P}^2}(0) =
- dim \ ker \ \Dd_{i,P} \ \ \text{for any $P \in
Gr_{\infty}^*(\Dd_i)$} \,,$$ hence we can use formula
(\ref{e:dt1}) in the situation we discuss under the assumption
$ker \ \Dd = \{0\}$.
\end{rem}

\vskip 1cm

We split $\z'_{\Dd_R^2}(0)$ into contributions coming from
different submanifolds plus cylinder contributions and the error
terms. Here $\Dd_R$ denotes the operator $\Dd$ on the manifold
$M_R$ equal to the manifold $M$ with $N$ replaced by $N_R$\,. We
introduce a manifold with boundary
$$M_{1,R} = M_1 \cup [-R,0] \times Y \,,$$
where we identify the ``old" collar neighborhood of the boundary
$Y$ on $M_1$ with $[-R-1,-R] \times Y$\,. Similarly we introduce
the manifold $M_{2,R}$\,. The bundle of Clifford modules $S$
splits on $Y$ into subbundles of spinors of positive and negative
chirality
$$S|_Y = S^+ \oplus S^- \ \ , \ \ \text{with  $S^{\pm} =
Ran \ {\frac 1{2}}(Id \mp i\G)$} \,.$$ The operator $P_{\pm} =
{\frac 1{2}}(Id \mp i\G)$ is the orthogonal projection of $S|_Y$
onto $S^{\pm}$ and provides $\Dd_i$ with a (local) {\it chiral}
elliptic boundary condition. This again means that the operator
$\Dd_{i,{\pm}} = \Dd_i$ with the domain
$$dom \ \Dd_{i,{\pm}} = \{s \in H^1(M_i;S) \mid P_{\pm}(s|_Y) = 0\}
\,,$$ is Fredholm and that its kernel and cokernel consist of only
smooth sections. We also have
\begin{equation}\label{e:dpm}
\Dd_{i,+}^* = \Dd_{i,-} \, .
\end{equation}
We study the $\z$-determinants of the corresponding Laplacians
\begin{equation}\label{e:pm1}
\Delta_{i,\pm} = \Dd_{i,\mp}\Dd_{i,\pm} \, .
\end{equation}
We denote by $\Delta_{i,R,\pm}$ the corresponding operator on the
manifold $M_{i,R}$\,.

\bigskip

In the present paper we avoid a discussion of the difficult issues
related to the existence of the {\it ``small"} eigenvalues of the
operators involved. Therefore we assume that the tangential
operator $B$ is invertible, i.e. $ker \ B = \{0\}\,.$ However,
this condition alone does not make all small eigenvalues
disappear. Careful analysis shows that we also need to assume that
the operator $\Dd_{i,\infty}$\,, equal to the operator $\Dd_i$
extended in a natural way to the manifold $M_{i,\infty}$\,, has no
$L^2$-solutions. The manifold $M_{i,\infty}$ is simply $M_i$ with
the infinite semicylinder $[0,\infty) \times Y$ (or $(-\infty,0]
\times Y$) attached (see \cite{CLM1}, see also \cite{KPW94}). The
existence of $L^2$-solutions of $\Dd_{i,\infty}$ on $M_{i,\infty}$
is responsible for the existence of exponentially small
eigenvalues of the operator $\Dd_R$. Therefore we assume
$ker_{L^2} \Dd_{i,\infty} = \{0\}$\,. The conditions we posed make
the small eigenvalues disappear. In particular, all the elliptic
boundary problems we discuss in this paper are invertible. We
refer to Proposition \ref{p:1} and  Remark \ref{r:lb} for more
information.

\bigskip

Our first main result is the following theorem

\bigskip

\begin{thm}\label{t:chiralsplit}
Let us assume that
\begin{equation}\label{e:ass1} ker_{L^2} \Dd_{1,\infty} = \{0\} =
ker_{L^2} \Dd_{2,\infty} \ \ \tand \ \ ker \ B = \{0\} \, .
\end{equation}

Then

\begin{equation}\label{e:chiralsplit1}
\lim_{R \to \infty} \{ln \ det_{\z}\Dd_R^2 - ln \
det_{\z}\Delta_{1,R,-} - ln \ det_{\z}\Delta_{1,R,+}\} = 0 \, ,
\end{equation}

\bigskip

or equivalently
\begin{equation}\label{e:chiralsplit2}
\lim_{R \to \infty}
{\frac{det_{\z}\Dd_R^2}{det_{\z}\Delta_{1,R,-}{\cdot}det_{\z}\Delta_{2
,R,+}}} = 1 \, .
\end{equation}
\end{thm}

\vskip 1cm

This Theorem is implicit in \cite{KlWo96}. The focus of the
authors was on the non-standard $\eta$-invariant introduced by
Singer in \cite{Si88} and on the analytic torsion. Therefore no
statement was made about the $\z$-determinant.

\bigskip

In Section 1 we use {\it Duhamel's Principle} to show that in
order to study
\begin{equation}\label{e:limit}
\lim_{R \to \infty} \{ln \ det_{\z}\Dd_R^2 - ln \
det_{\z}\Delta_{1,R,+} - ln \ det_{\z}\Delta_{2,R,+}\} \,
\end{equation}

\noi it is enough to discuss the cylinder contributions.

In Section 2 we perform the computation on the cylinder and show
that the limit (\ref{e:limit}) is indeed equal to $0$ . Then we
study the difference between the cylinder contribution for the
chiral boundary condition and for the Atiyah-Patodi-Singer
condition. Straightforward computations show that a new term
appears which is equal to $-ln \ 2{\cdot}\z_{B^2}(0)$\,. This
gives the main result of the paper:

\bigskip

\begin{thm}\label{t:chiralsplit2}
The following equality holds under the assumptions of our Theorem
\ref{t:chiralsplit}
\begin{equation}\label{e:chiralsplit3}
\lim_{R \to \infty}
{\frac{det_{\z}\Dd_R^2}{det_{\z}\Dd_{1,R,\Pi_<}^2
{\cdot}det_{\z}\Dd_{2,R,\Pi_>}^2}} = 2^{-\z_{B^2}(0)} \, .
\end{equation}

\end{thm}

\bigskip
The Appendix by Yoonweon Lee contains a refined version of the
computations of the cylinder contribution to the trace of the heat
kernel of the Atiyah-Patodi-Singer problem performed by the second
author in \cite{KPW99}. The more careful analysis by Lee proves
$mod \ \Z$ vanishing of the function \hskip 2mm $P \mapsto
\z_{\Dd_P^2}(0)$ on the Grassmannian $Gr_{\infty}^*(\Dd_i)$ .
Moreover, the formula (\ref{e:fo2}) (see Appendix Proposition
\ref{p:fo2}) is used in the proof of Theorem \ref{t:chiralsplit2}.

\bigskip

\begin{rem}\label{r:inv1}
This paper is related to many other works on the gluing formulas
for the $\z$-determinants. We refer to an excellent survey article
\cite{MP98} for the review of different approaches and the
extensive bibliography. However, we want to mention that Theorem
\ref{t:chiralsplit2} is closely related to the results of
\cite{Ha98}. In \cite{Ha98} only the operator $d + d^*$ is
treated, but the gluing formula similar to (\ref{e:chiralsplit3})
is obtained using the b-calculus technique, in the situation where
the zero eigenvalues are allowed.
\end{rem}

\bigskip

\section{Duhamel's Principle. Reduction to
the Cylinder} \label{s:bl}

\vskip 5mm

Our assumptions about the operator $\Dd_R$ (see (\ref{e:ass1}))
allow us to apply the technique developed in \cite{DoWo91} and to
reduce the proof of Theorem \ref{t:chiralsplit} and Theorem
\ref{t:chiralsplit2} to the computations on the cylinder. The
first important Corollary of (\ref{e:ass1}) is the following
Proposition

\vskip 5mm

\begin{prop}\label{p:1}
Let us assume that (\ref{e:ass1}) holds. Then there exist positive
constants $c$ and $R_0$\,, such that
\begin{equation}\label{e:lb}
\mu > c
\end{equation}
for any eigenvalue $\mu$ of the operator $\Dd_R^2$\,,
$\Delta_{i,R,\pm}$\,, $\Dd_{1,R,\Pi_<}^2$\,, $\Dd_{2,R,\Pi_>}^2$\
and for any $R > R_0$\,.
\end{prop}

\bigskip

\begin{rem}\label{r:lb}
The estimate (\ref{e:lb}) was observed by W. M${\rm {\ddot
u}}$ller. We refer to \cite{DoWo91} (Theorem 6.1) for the proof in
the case of the Atiyah-Patodi-Singer condition (operators
$\Dd_{1,R,\Pi_<}^2$, $\Dd_{2,R,\Pi_>}^2$). A more general result
was published in \cite{Mu94}, Proposition 8.14. The proof for the
``chiral" boundary conditions (operators $\Delta_{i,R,\pm}$) is
even more simple. The case of the operator $\Dd_R$ was analyzed in
\cite{CLM1} (see also \cite{KPW94}).
\end{rem}

\vskip 5mm

We need to recall the following result

\vskip 5mm

\begin{prop}\label{p:2}
Let $\Ee_R(t;x,y)$ denote the kernel of the {\it heat operator}
for $\Delta_R$\,, where $\Delta_R$ denotes one of the operators
from Proposition \ref{p:1}. Assume that (\ref{e:ass1}) holds. Then
there exist positive constants $c_1$ , $c_2$ and $c_3$ such that
\begin{equation}\label{e:lb1}
\|\Ee_R(t;x,y)\| \le c_1t^{-{\frac
n{2}}}e^{c_2t}e^{-c_3{\frac{d^2(x,y)}{t}}} \, ,
\end{equation}
for any $t > 0$ and any $x,y \in M_R$ ($M_{1,R}$\,, or $M_{2,R}$
respectively) and for any $R > R_0$\,.
\end{prop}

\vskip 5mm

We refer to Sections 2 and 4 of \cite{DoWo91} for the proof and
related results. In particular Proposition \ref{p:2} implies the
following estimate
\begin{equation}\label{e:volume}
|Tr \ e^{-\Delta_R}| < c_4{\cdot}R \, .
\end{equation}

\bigskip

Now, we are ready to prove that we can neglect the {\it ``large
time contribution"} to the $\z$-determinant of $\Delta_R$\,.

\bigskip

\begin{cor}\label{c:lb1}
Let us assume (\ref{e:ass1})\,, then for any $\e >0$ the following
equality holds
\begin{equation}\label{e:ltc}
\lim_{R \to \infty} \int_{R^{\e}}^{\infty}{\frac 1{t}}{\cdot}Tr \
e^{-t\Delta_R}dt = 0 \, .
\end{equation}
\end{cor}

\bigskip

\begin{proof}
Assume that $R > R_0$ and let $\{\mu_k\}_{k=1}^{\infty}$ denote
the set of eigenvalues of $\Delta_R$\,. We have

\begin{align*}
\int_{R^{\e}}^{\infty}{\frac 1{t}}{\cdot}Tr \ e^{-t\Delta_R}dt &=
\int_{R^{\e}}^{\infty}{\frac 1{t}}{\cdot} \sum_{k =
1}^{\infty}e^{-t\mu_k}dt = \int_{R^{\e}}^{\infty}{\frac
1{t}}{\cdot}\sum_{k =
1}^{\infty}e^{-(t-1)\mu_k}e^{-\mu_k}dt\\
& < \int_{R^{\e}}^{\infty}{\frac 1{t}}e^{-(t-1)c}{\cdot}Tr \
e^{-\Delta_R}dt \ \ ,
\end{align*}
where $c$ is the constant from Proposition \ref{p:1}. We use
(\ref{e:volume})
$$\int_{R^{\e}}^{\infty}{\frac 1{t}}{\cdot}Tr \
e^{-t\Delta_R}dt < \int_{R^{\e}}^{\infty}{\frac
1{t}}e^{-(t-1)c}{\cdot}Tr \ e^{-\Delta_R}dt <
c_6R^{1-\e}{\cdot}e^{-c_7R^{\e}}$$ and the Corollary follows
easily.
\end{proof}

\vskip 5mm

Now we follow \cite{DoWo91}. Let ${\tilde \Ee}_R(t;x,y)$ denote
the kernel of the operator $e^{-t\Dd_R^2}$ on the manifold $M_R$
and let $\Ee_{cyl}(t;x,y)$ denote the kernel of the operator
$e^{-t(-\partial_u^2 + B^2)}$ on the  infinite cylinder
$(-\infty\,, +\infty) \times Y$, or the kernel of the APS -
operator on $((-\infty, 0] \cup [0, \infty)) \times Y$. We
introduce a smooth, increasing function $\rho(a,b) : [0, \infty)
\to [0,1]$ equal to $0$ for $0 \le u \le a$ and equal to $1$ for
$b \le u$\,. We use $\rho(a,b)(u)$ to define
$$\phi_1 = 1 - \rho({\frac 5{7}}R\,, {\frac
6{7}}R) \ \ , \ \ \psi_1 = 1 - \psi_2$$ and
$$\phi_2 = \rho({\frac
1{7}}R\,, {\frac 2{7}}R) \ \ , \ \ \psi_2 = \rho({\frac 3{7}}R\,,
{\frac 4{7}}R) \,.$$ We extend these functions to symmetric
functions on the whole real line. All these functions are constant
outside the interval $[-R\,, R]$ and we use them to define the
corresponding functions on the manifold $M_R$\,. Now we define
$Q_R(t;x,y)$ as a {\it ``parametrix"} for the kernel ${\tilde
\Ee}_R(t;x,y)$ \,, actually using ${\tilde \Ee}_R(t;x,y)$\,, but
the point here is that we will be able to separate the cylinder
and the interior contribution
\begin{equation}\label{e:p1}
Q_R(t;x,y) = \phi_1(x){\Ee}_{cyl}(t;x,y)\psi_1(y) +
\phi_2(x){\tilde \Ee}_R(t;x,y)\psi_2(y) \, .
\end{equation}
Standard computations show that
\begin{equation}\label{e:p2}
{\tilde \Ee}_R(t;x,y) = Q_R(t;x,y) + ({\tilde \Ee}_R*\Cc_R)(t;x,y)
\, ,
\end{equation}

\noi where ${\tilde \Ee}_R*\Cc_R$ is a convolution given by
$$({\tilde \Ee}_R*\Cc_R)(t;x,y) = \int_0^tds\int_{M_R}dz \
{\tilde \Ee}_R(s;x,z)\Cc_R(t-s;z,y) \,,$$ and the correction term
$\Cc_R(t;x,y)$ is given by the formula
\begin{multline}\label{e:p3}
\Cc_R(t;x,y) =
-{\frac{{\partial}^2\phi_1}{{\partial}u^2}}(x){\Ee}_{cyl}(t;x,y)\psi_1(y)
- {\frac{{\partial}\phi_1}{{\partial}u}}(x)
{\frac{{\partial}{\Ee}_{cyl}}{{\partial}u}}(t;x,y)\psi_1(y)
\\
-{\frac{{\partial}^2\phi_2}{{\partial}u^2}}(x){\tilde
\Ee}_R(t;x,y)\psi_2(y) - {\frac{{\partial}\phi_2}{{\partial}u}}(x)
{\frac{{\partial}{\tilde \Ee}_R}{{\partial}u}}(t;x,y)\psi_2(y) \,.
\end{multline}

\bigskip
\noi The choice of the cut-off functions and the estimate
(\ref{e:lb1}) allow us to neglect the {\it ``error" term
contribution} to the logarithm of the determinant in the limit as
$R \to \infty$.

\bigskip

\begin{lem}\label{l:err1}
The error term $\Cc_R(t;x,y)$ is equal to $0$ outside of the
cylinder $[-{\frac 6{7}}R\,, {\frac 6{7}}R] \times Y$\,, moreover
it is equal to $0$ if the distance between $x$ and $y$ is smaller
than ${\frac R{7}}$\,. This fact combined with (\ref{e:lb1})
proves the following estimate for certain positive constants
\begin{equation}\label{e:err1}
\|({\tilde \Ee}_R*\Cc_R)(t;x,y) \| \le
c_1e^{c_2t}e^{-c_3{{\frac{R^2}{t}}}} \, .
\end{equation}
\end{lem}

\bigskip

The proof is elementary and follows the proof of the similar
statement in \cite{DoWo91} (see Proposition 5.2 of \cite{DoWo91}).

\bigskip

\begin{cor}\label{c:err1}
Assume that $0 < \e < 1$\,, then
\begin{equation}\label{e:err2}
\lim_{R \to 0}\int_0^{R^{\e}}{\frac{dt}{t}}\int_{\Mm_R}tr \
({\tilde \Ee}_R*\Cc_R)(t;x,x)dx = 0 \, ,
\end{equation}
where $\Mm_R$ denotes any of the manifolds on which the operator
$\Delta_R$ (of Proposition 1.3) is defined.
\end{cor}

\bigskip

\begin{proof}
This follows from the estimate on the kernel $({\tilde
\Ee}_R*\Cc_R)(t;x,x)$
\begin{align*}
|tr \ ({\tilde \Ee}_R*\Cc_R)(t;x,x)| &\le \|({\tilde
\Ee}_R*\Cc_R)(t;x,x)\|\\
&\le \int_0^tds\int_{\Mm_R}\|{\tilde
\Ee}_R(s;x,z)\Cc_R(t-s;z,x)\|dz\\
&\le c_1e^{c_2t}{\cdot}\int_0^tds
\int_{\Mm_R}e^{-c_3{\frac{d^2(x,z)}{s}}}
e^{-c_3{\frac{d^2(x,z)}{t-s}}}dz \,.
\end{align*}
We use Lemma \ref{l:err1}. It follows that the integral with
respect to $z$ is taken over the cylinder and moreover that the
distance $d(x,z)$ is always larger than ${\frac R{7}}$\,, which
gives
\begin{multline*}
|tr \ ({\tilde \Ee}_R*\Cc_R)(t;x,x)| \le
c_1e^{c_2t}{\cdot}\int_0^tds
\int_{\Mm_R}e^{-c_3{\frac{d^2(x,z)}{s}}}e^{-c_3{\frac{d^2(x,z)}{t-s}}}dz\\
< c_1e^{c_2t}\int_0^tds\int_{-R}^Re^{-c_4{\frac{tR^2}{s(t-s)}}}dz
< c_5Re^{c_2t}e^{-c_6{\frac{R^2}{t}}}\int_0^tds <
c_5Rte^{c_2t}e^{-c_6{\frac{R^2}{t}}} \,.
\end{multline*}
Now we have
\begin{multline*}
\Bigl|\int_0^{R^{\e}}{\frac{dt}{t}}\int_{\Mm_R}tr \ ({\tilde
\Ee}_R*\Cc_R)(t;x,x)dx\Bigr|\\ <
\int_0^{R^{\e}}{\frac{dt}{t}}c_5Rte^{c_2t}e^{-c_6{\frac{R^2}{t}}}
< c_5R^{1+\e}e^{-c_7R^{1 - \e}} \,,
\end{multline*}
and (\ref{e:err2}) is proved.
\end{proof}

\vskip 5mm

The last result clearly explains that we have to analyze only the
cylinder contribution to study the ratio of the determinants in
the adiabatic limit. This is done in the next Section.

\vskip 5mm

\section{Computations on the Cylinder}
\label{s:cl}

\vskip 5mm

Our study of the decomposition formula for the $\z$-determinant
involves the {\it ``Laplacians"} $\Delta_{i,\pm} =
\Dd_{i,\mp}\Dd_{i,\pm}$\,.  It is well-known that $\Delta_{i,+}$
is the operator $\Dd_i^2$ subject to the Dirichlet boundary
condition on $S^+$ and the Neumann boundary condition on $S^-$
(see for instance Lemma 1.1. in \cite{KlWo96}). %%
%% \vfill\eject

It was explained in the previous section that it is enough to
analyze the cylinder contribution. Hence we have to write down the
kernel of the {\it heat operator} defined by
$$-\partial_u^2 + B^2 : C^{\infty}([0\,, \infty) \times Y ; S = S^+
\oplus S^-) \to  C^{\infty}([0\,, \infty) \times Y ; S = S^+
\oplus S^-)$$

\noi subject to the Dirichlet condition on ${S^+}|_{\{0\} \times
Y}$ and the Neumann condition on ${S^-}|_{\{0\} \times Y}$\,, in
the case of the operator coming from the manifold $M_{2,R}$\,; and
subject to the Dirichlet condition on ${S^-}|_{\{0\} \times Y}$
and the Neumann condition on ${S^+}|_{\{0\} \times Y}$\,, in the
case of the operator  coming from the manifold $M_{1,R}$\,.

\noi Let $\Ee^+(t;x,y)$ denote the kernel of the first operator.
The explicit formula is well-known (see \cite{KlWo96} for
instance)
\begin{align}\label{e:+heat}
\Ee^+(t;(u,x),(v,y)) =&
{\frac1{\sqrt{4{\pi}t}}}\{e^{-\frac{(u-v)^2}{4t}} -
e^{-\frac{(u+v)^2}{4t}}\}
e^{-tB^2}(t;x,y)P_+(y) \\
&+ {\frac1{\sqrt{4{\pi}t}}}\{e^{-\frac{(u-v)^2}{4t}} +
e^{-\frac{(u+v)^2}{4t}}\} e^{-tB^2}(t;x,y)P_-(y) \,, \notag
\end{align}

\bigskip
\noi where $e^{-tB^2}(t;x,y)$ denotes the kernel of the operator
$e^{-tB^2}$\,. This formula determines the cylinder contribution
coming from the manifold $M_{2,R}$\,. The inward normal coordinate
on $M_1$ is equal to $-u$. As a consequence the chirality of
spinors, which is determined by the Clifford multiplication by the
normal vector, is switched as $G$ is replaced by $-G$\,. The
corresponding cylinder contribution for the manifold $M_{1,R}$  is
determined by the kernel
\begin{align}\label{e:-heat}
\Ee^-(t;(u,x),(v,y)) =& {\frac1{\sqrt{4{\pi}t}}}\{e^{-
\frac{(u-v)^2}{4t}} + e^{- \frac{(u+v)^2}{4t}}\}
e^{-tB^2}(t;x,y)P_+(y) \\
&+ {\frac1{\sqrt{4{\pi}t}}}\{e^{- \frac{(u-v)^2}{4t}} -
e^{-\frac{(u+v)^2}{4t}}\} e^{-tB^2}(t;x,y)P_-(y) \,. \notag
\end{align}

\bigskip

\noi Now we sum up the formulas (\ref{e:+heat}) and
(\ref{e:-heat}) and put $u=v$\,, $x=y$\,. This gives

\begin{multline}\label{e:heat2}
\int_0^{R}du{\frac 1{\sqrt{4{\pi}t}}}\{1 -
e^{-{\frac{u^2}{t}}}\}Tr_Ye^{-tB^2}P_+ +
\int_0^{R}du{\frac1{\sqrt{4{\pi}t}}}\{1 +
e^{-{\frac{u^2}{t}}}\}Tr_Ye^{-tB^2}P_-
+\\
\int_0^{R}du{\frac1{\sqrt{4{\pi}t}}}\{1 +
e^{-{\frac{u^2}{t}}}\}Tr_Ye^{-tB^2}P_+ +
\int_0^{R}du{\frac1{\sqrt{4{\pi}t}}}\{1 -
e^{-{\frac{u^2}{t}}}\}Tr_Ye^{-tB^2}P_- \,.
\end{multline}

\bigskip
\noi In the formula (\ref{e:heat2}) we neglect the presence of the
cut-off functions introduced in Section 1.  We also denote by
$Tr_Ye^{-tB^2}$ the trace of the operator $e^{-tB^2}$ on the
manifold $Y$\,. Therefore modulo a term exponentially decaying
with $R$, the boundary contribution to the sum \hskip 2mm $Tr \
e^{-t\Delta_{R,1}} + Tr \ e^{-t\Delta_{R,2}}$ is equal to

{\begin{equation}\label{e:bdctr} {\frac
2{\sqrt{4{\pi}t}}}{\cdot}\int_0^{R}du\,(Tr_Ye^{-tB^2}P_+ +
Tr_Ye^{-tB^2}P_-) = {\frac
1{\sqrt{4{\pi}t}}}{\cdot}\int_{-R}^{R}du{\cdot}Tr_Ye^{-tB^2} \, .
\end{equation}

\bigskip

\noi The right side of (\ref{e:bdctr}) is exactly equal to the
trace of the heat kernel of the operator $-\partial_u^2 + B^2$ on
the cylinder $(-\infty, +\infty) \times Y$\,, which is the
cylinder contribution of the operator $\Dd_R^2$ modulo terms which
disappear as $R \to \infty$. This ends the proof of Theorem
\ref{t:chiralsplit}.

\vskip 5mm

Now, we have to analyze the difference between the trace of
$\Ee^+(t;x,y)$ and the trace $\Ee^>(t;x,y)$\,, where
$\Ee^>(t;x,y)$ denotes the kernel of the {\it heat operator}
defined by the operator $G(\partial_u + B)$ subject to the
Atiyah-Patodi-Singer boundary condition. We introduce $\phi(u)$ a
smooth cut-off function, equal to $1$ for $0 \le u \le R$ and
vanishing for $2R \le u$\,, with derivatives bounded by ${\frac
c{R}}$\,, and we study the following function
\begin{multline}\label{e:ar1}
\Tt(s)=\int_0^{\infty}t^{s-1}dt\int_{ [0\,, \infty) \times
Y}\phi(u){\cdot}tr(\Ee^>(t;(u,y),(u,y))\\
- \Ee^+(t;(u,y),(u,y)))dydu \,.
\end{multline}
Long, but elementary computations give us the following formula
for the contribution made by the Atiyah-Patodi-Singer part (see
Appendix, Proposition A.4.)
\begin{align}\label{e:ar2}
&\int_0^{\infty}t^{s-1}dt\int_{[0\,, \infty) \times
Y}\phi(u){\cdot}tr(\Ee^>(t;(u,y),(u,y))dydu\\
     =& {\frac
1{\sqrt{4\pi}}}\int_0^{\infty}\phi(u)du\int_0^{\infty}t^{s -
{\frac
3{2}}}Tr_Ye^{-tB^2}dt \notag \\
    & + {\frac
1{2}}{\cdot}\int_0^{\infty}t^{s-1}dt
\int_0^{\infty}{\phi}'(u)du\sum_{n = 1}^{\infty} e^{2u{\la}_n}{\rm
erfc}\Bigl({\frac u{\sqrt{t}}} + {\la}_n\sqrt{t}\Bigr)\notag
\\
&+ {\frac{\G(s + {\frac1{2}})}{4s\sqrt{\pi}}}\z_{B^2}(s) - {\frac
1{2\sqrt{\pi}}}{\cdot}\int_0^{\infty}t^{s - {\frac 3{2}}}\,Tr_Y
e^{-tB^2}
\Bigl(\int_0^{\infty}\phi(u)e^{-{\frac{u^2}{t}}}du\Bigr)\,dt \,.
\notag
\end{align}

\bigskip
\noi We have three terms on the right side of (\ref{e:ar2}), which
we denote by $\Tt_1(s)$, $\Tt_2(s)$ and $\Tt_3(s)$\,. The sum in
$\Tt_2(s)$ is taken over all positive eigenvalues $\la_n$ of the
tangential  operator $B$ and the function ${\rm erfc}(u)$ is given
by the formula
\begin{equation}\label{e:erfc}
{\rm erfc}(u) = {\frac 2{\sqrt{\pi}}}\int_u^{\infty}e^{-s^2}ds \,.
\end{equation}
The first contribution $\Tt_1(s)$ corresponds to the contribution
to (\ref{e:ar1}) given by the kernel $\Ee^+(t;(u,y),(u,y))$ and
they cancel each other when we take the difference. We can also
easily deal with the second contribution:

\bigskip

\begin{prop}\label{p:ar1}
The function
$$\Tt_2(s) = {\frac 1{2}}
{\cdot}\int_0^{\infty}t^{s-1}dt\int_0^{\infty}{\phi}'(u)du
\sum_{n=1}^{\infty}e^{2u{\la}_n}{\rm erfc}({\frac u{\sqrt{t}}} +
{\la}_n\sqrt{t})$$

\bigskip
\noi is a holomorphic function of $s$ vanishing as $R \to
\infty$\,.
\end{prop}

\vskip 5mm

\begin{proof}
We estimate using $\int_u^{\infty}e^{-s^2}ds \le e^{-u^2}$

$$\int_0^{\infty}t^{s-1}dt\int_0^{\infty}{\phi}'(u)du\sum_{n
= 1}^{\infty}e^{2u{\la}_n}{\rm erfc}({\frac u{\sqrt{t}}} +
{\la}_n\sqrt{t}) \le$$

$${\frac
2{\sqrt{\pi}}}\int_0^{\infty}t^{s-1}dt\int_R^{2R}{\phi}'(u)du\sum_{n
= 1}^{\infty}e^{-{\frac{u^2}{t}}}e^{-t\la_n^2} \le$$

$${\frac c{R}}{\frac
2{\sqrt{\pi}}}\int_0^{\infty}t^{s-1}(\sum_{n =
1}^{\infty}e^{-t\la_n^2})dt\int_R^{2R}e^{-{\frac{u^2}{t}}}du \le
$$
$$
{\frac{c_1}{R}}\int_0^{\infty}t^{s-1}Tr \
e^{-tB^2}dt\int_R^{\infty}e^{-{\frac{u^2}{t}}}du \le$$

$${\frac{c_1}{R}}\int_0^{\infty}t^{s-{\frac1{2}}}Tr
\
e^{-tB^2}dt\int_R^{\infty}e^{-{\frac{u^2}{t}}}{\frac{du}{\sqrt{t}}}
\le {\frac{c_1}{R}}\int_0^{\infty}t^{s-{\frac 1{2}}}
e^{-{\frac{R^2}{t}}}Tr \ e^{-tB^2}dt \,,$$

\bigskip
\noi and the Proposition follows.
\end{proof}

\vskip 5mm Now we see that $\Tt_3(s)$ is the only source of an
additional contribution. It is not difficult to see that, modulo a
function holomorphic on the whole complex plane, $\Tt_3(s)$ is
equal to
\begin{align*}
\Ss(s) &= {\frac{\G(s + {\frac1{2}})}{4s\sqrt{\pi}}}\z_{B^2}(s) -
{\frac1{4}}{\cdot}\int_0^{\infty}t^{s - 1} Tr_Y e^{-tB^2}dt\\
&= {\frac{\G(s + {\frac 1{2}})}{4s\sqrt{\pi}}}\z_{B^2}(s) -
{\frac{\G(s)}{4}}{\cdot}\z_{B^2}(s) \,.
\end{align*}
Indeed, the difference
$$g_R(s) = \Tt_3(s) - \Ss(s) =
{\frac 1{2\sqrt{\pi}}}\int_0^{\infty}t^{s-1}Tr \
e^{-tB^2}(\int_0^{\infty}(1-\phi(u))
e^{-{\frac{u^2}{t}}}{\frac{du}{\sqrt{ t}}})dt$$ is a holomorphic
function on the complex plane, which depends on the parameter
$R$\,. We use the following elementary result:

\bigskip

\begin{lem}\label{l:ar1}
The following equality holds for any complex $s$
\begin{equation}\label{e:ar3}
\lim_{R \to \infty} g_R(s) = \lim_{R \to \infty} g_R'(s) = 0 \, .
\end{equation}
\end{lem}

\bigskip

\begin{proof}
We have to estimate
\begin{equation}\label{e:a}
\Bigl|{\frac1{2\sqrt{\pi}}}\int_0^{\infty}t^{s-1}Tr \
e^{-tB^2}\int_0^{\infty}\bigl(1-\phi(u)\bigr)
e^{-{\frac{u^2}{t}}}{\frac{d u}{\sqrt{t}}}dt\Bigr| \, .
\end{equation}
We use the following elementary inequality

\begin{align*}
\Bigl|\int_0^{\infty}(1-\phi(u))e^{-{\frac{u^2}{t
}}}{\frac{du}{\sqrt{t}}}\Bigr| &\le
\Bigl|\int_R^{\infty}e^{-{\frac{u^2}{t}}
}{\frac{du}{\sqrt{t}}}\Bigr| =
\Bigl|\int_{\frac{R}{\sqrt{t}}}^{\infty}(-{\frac{1}{2s}})(-2s)e^{-s^2}
ds\Bigr|\\
&\le
\Bigl|-{\frac{\sqrt{t}}{2R}}\int_{\frac{R}{\sqrt{t}}}^{\infty}{
\frac{d}{ds}} (e^{-s^2})ds\Bigr| =
{\frac{\sqrt{t}}{2R}}e^{-{\frac{R^2}{t}}} \,.
\end{align*}
This allows us to estimate (\ref{e:a})
\begin{multline*}
\Bigl|{\frac{1}{2\sqrt{\pi}}}\int_0^{\infty}t^ {s-1}Tr \
e^{-tB^2}\int_0^{\infty}(1-\phi(u))e^{-{\frac{u^2}{t}}}{\frac{du}{\sqrt{t}
}}
dt\Bigr|\\
< {\frac
1{4{\sqrt{\pi}}R}}\int_0^{\infty}t^{s-{\frac1{2}}}e^{-{\frac{R^2}{t}}}
Tr_Y e^{-tB^2}dt \,.
\end{multline*}
%% \vfill\eject
%%
%% \noi
The last expression goes to $0$ as $R \to \infty$\,. The estimates
on the derivatives with respect to $s$ go exactly in the same way.
\end{proof}

\vskip 5mm

The function $\Ss(s)$ was given by the formula
$$\Ss(s) = \Bigl({\frac{\G(s + {\frac1{2}})}{4s\sqrt{\pi}}}
- {\frac{\G(s)}{4}}\Bigr){\cdot}\z_{B^2}(s) \,.$$ We see that
$\Ss(s)$ is a holomorphic function for $Re(s) > {\frac n{2}}$ and
that it has a meromorphic extension to the whole complex plane
with simple poles on the real line, provided by both factors.
Hence the poles at the positive integers come from $\z_{B^2}(s)$
and the $\z$-function is regular in the neighborhood of $0$\,. The
first factor
$${\frac{\G(s +
{\frac1{2}})}{4s\sqrt{\pi}}} - {\frac{\G(s)}{4}}$$ is holomorphic
for $Re(s) > 0$ and it is not very difficult to show that in fact
it is holomorphic in a neighborhood of $s=0$\,. We have
\begin{align*}
{\frac{\G(s + {\frac1{2}})}{4s\sqrt{\pi}}} - {\frac{\G(s)}{4}} &=
{\frac1{4\sqrt{\pi}}}{\cdot} {\frac{\G(s +
{\frac1{2}}) - s\G(s)\sqrt{\pi}}{s}} \\
&= {\frac 1{4\sqrt{\pi}}}{\cdot}\Bigl({\frac{\G(s + {1\over{2}}) -
\G({1\over{2}})}{s}} + \G({1\over{2}}){\frac{1 -
\G(s+1)}{s}}\Bigr) \,,
\end{align*}
and we see that
$$\lim_{s \to 0}{\frac{\G(s +
{1\over{2}})}{4s\sqrt{\pi}}} - {\frac{\G(s)}{4}} = {\frac
1{4\sqrt{\pi}}}{\cdot} \Bigl(\G'({\frac 1{2}}) -
\sqrt{\pi}{\cdot}\G'(1)\Bigr) \,.$$ It is well-known that $\G'(1)
= \g$ (once again, $\g$ denotes Euler's constant), and it is not
difficult to compute $\G'({1\over{2}})$ using, for instance, the
formula
$$\G(z + {\frac1{2}})
= {\frac{{\sqrt{\pi}}{\cdot}\G(2z)}{2^{2z-1}{\cdot}\G(z)}} \,,$$
(see for instance \cite{Ta196}, formula (A22) on page 265).

\bigskip
\noi We have
\begin{align*}
&\lim_{s \to 0}{\frac{\G(s + {\frac 1{2}}) - \G({\frac1{2}})}{s}}
= \sqrt{\pi}{\cdot} \lim_{s \to 0}
{\frac{{{\G(2s)}/{2^{2s-1}{\cdot}\G(s)}} - 1}{s}}\\
=& \sqrt{\pi}{\cdot}\lim_{s \to 0} {\frac{2^{1-2s}\G(2s) -
\G(s)}{s\G(s)}} = \sqrt{\pi}{\cdot}\lim_{s \to 0}{\frac
1{\G(1+s)}}{\cdot} \lim_{s \to
0}(2^{1-2s}\G(2s) - \G(s)) \\
=& \sqrt{\pi}{\cdot}\lim_{s \to 0}(2^{1-2s}({\frac1{2s}} + \g +
2sh(2s)) - ({\frac 1{s}} + \g + sh(s))) \,,
\end{align*}
where $h(s)$ is a holomorphic function in the neighborhood of
$s=0$\,. Hence we finally obtain
\begin{align*}
&\lim_{s \to
0}{\frac{\G(s + {\frac 1{2}}) - \G({\frac 1{2}})}{s}}\\
     =&
\sqrt{\pi}{\cdot}\lim_{s \to 0} \Bigl({\frac{2^{1-2s} - 2}{2s}} +
2^{1-2s}\g - \g + 2^{1-2s}2s\,h(2s) -
s\,h(s)\Bigr)\\
=&-2\sqrt{\pi}{\cdot}ln \ 2 +\sqrt{\pi}\gamma \,,
\end{align*}
and
\begin{equation}\label{e:r1}
\lim_{s \to 0}{\frac{\G(s + {\frac1{2}})}{4s\sqrt{\pi}}} -
{\frac{\G(s)}{4}} = -
{\frac1{4\sqrt{\pi}}}{\cdot}2\sqrt{\pi}{\cdot}ln \ 2 = -
{\frac1{2}}ln \ 2 \, .
\end{equation}

\noi This gives us the following result

\bigskip

\begin{prop}\label{p:r1}
The adiabatic limit of the difference between the logarithm of the
$\z$-determinant of the operator $\Dd_{2,R,\Pi_>}^2$ and the
logarithm of the $\z$-determinant of the operator $\Delta_{2,R,+}$
is given by
\begin{equation}\label{e:r2}
\lim_{R \to \infty} (ln \ det_{\z}\Dd_{2,R,\Pi_>}^2 - ln \
det_{\z}\Delta_{2,R,+}) = {\frac{ln \ 2}{2}}{\cdot}\z_{B^2}(0) \,
.
\end{equation}
\end{prop}

\vskip 1cm We have obtained {\it ``half"} of the correction term
which appears in Theorem \ref{t:chiralsplit2} (see
(\ref{e:chiralsplit3})). The other {\it ``half"} is equal to the
contribution of the manifold $M_{1,R}$\,. Now Theorem
\ref{t:chiralsplit2} is proved.

\vskip 1 cm

\begin{appendix}
\section{The value
of the $\zeta$-function at $s=0$ \\on the smooth, self-adjoint
Grassmannian}

\begin{center}
{\sc Yoonweon Lee}
\end{center}

\vskip 5mm {\bf Acknowledgements.} The author was supported by
Korea Research Foundation Grant KRF-2000-015-DP0045.

\vskip 5mm

In this Appendix we write $M$ instead of $M_2$ and $\Dd$ instead
of $\Dd_2$\,. The goal is to prove the following result

\bigskip

\begin{prop}\label{p:zeta}
For any $P \in Gr_{\infty}^*(\Dd)$, the value of
$\z_{\Dd_{P}^2}(s)$ at $s = 0$ is equal to $- dim \ ker \Dd_P$ \,.
\end{prop}

\bigskip

\begin{rem}\label{r:zeta}

(1) The proof depends only on the assumption that the perturbation
of $\Pi_>$ is an operator of the trace class. Therefore the result
holds for any orthogonal projection $P = Id + GPG$\,, such that $P
- \Pi_>$ is a pseudodifferential operator of order $- dim \ Y -
1$\,.

(2) One of the formulas we obtain for the cylinder contribution to
the invariant $\z_{\Dd_{\Pi_>}^2}(0)$ (see Proposition A.4.) is
used in the proof of the decomposition formula for the
$\z$-determinant discussed in the main body of the paper.
\end{rem}

\bigskip

We start with the proof of Proposition \ref{p:zeta} in the most
simple case. We assume

\begin{equation}\label{e:zz1}
dim \ ker \ B = 0 \ \ and \ \ dim \ ker \ \Dd_{\Pi_>} = 0 \ \, .
\end{equation}

\bigskip

It was explained earlier  that the first condition in
(\ref{e:zz1}) implies that $\Pi_> \in Gr_{\infty}^*(\Dd)$ , hence
$\Dd_{\Pi_>}$ is a self-adjoint operator. The second condition
implies the invertibility of $\Dd_{\Pi_>}$ . We have to show that
$\z_{\Dd_{\Pi_>}^2}(0) = 0$ .

We start with selecting a smooth cut-off function $\rho : M \to
[0,1]$ equal to $1$ on $[0\,, {1\over{3}}] \times Y$ and equal to
$0$ on $M \setminus ([0\,,{2\over{3}}] \times Y)$\,. We also
choose $\rho_1\,, \rho_2 : M \to [0,1]$ such that
$$\rho_1|_{supp \ \rho} \equiv 1 \ \
and \ \ \rho_1 \equiv 0 \ on \ {M \setminus N} \ \ and$$
$$\rho_2|_{supp \ (1 -\rho)} \equiv 1 \ \ and \ \
\rho_2 \equiv 0 \ on \ {[0,{1\over{4}}] \times Y}  \,.$$ Let
$\Ee_{cyl}(t;x,y)$ denote the heat kernel of the
Atiyah-Patodi-Singer problem on the cylinder $[0\,, \infty) \times
Y$ and ${\tilde \Ee}(t;x,y)$ denote the kernel of the operator
$e^{-t{\tilde \Dd}^2}$\,, where $\tilde \Dd$ is the double of the
operator $\Dd$, living on $\tilde M$ the double of $M$ (see
\cite{BoWo93} for details of the construction). Finally let
$\Ee_>(t;x,y)$ denote the kernel of the {\it heat operator} of
$\Dd_{\Pi_>}^2$ on $M$\,. A standard application of {\it Duhamel's
Principle} shows that there exists a positive constant $c$\,, such
that

\begin{equation}\label{e:z1}
\Ee_>(t;x,y) = \rho_1(x)\Ee_{cyl}(t;x,y)\rho(y) + \rho_2(x){\tilde
\Ee}(t;x,y)(1 - \rho(y)) + O(e^{-{c\over{t}}}) \, ,
\end{equation}

\bigskip

\noi for $0 < t \le 1$\,. Now the $\z$-function is given by the
formula
$$\z_{\Dd_{\Pi_>}^2}(s) =
{\frac{1}{\G(s)}}\int_0^{\infty}t^{s-1}Tr \ e^{-t\Dd_{\Pi_>}^2}dt
= {\frac{1}{\G(s)}}\int_0^{\infty}t^{s-1}dt\int_M tr \
\Ee_>(t;x,x)dx \,.$$

\bigskip
\noi Equation (\ref{e:z1}) implies that there exist positive
constants $c_1$ and $c_2$\,, such that
\begin{equation}\label{e:z2}
|tr \ \Ee_>(t;x,x) - \rho(x){\cdot}tr \ \Ee_{cyl}(t;x,x) - (1 -
\rho(x)){\cdot}tr \ {\tilde \Ee}(t;x,x)| < c_1e^{-{{c_2}\over{t}}}
\,
\end{equation}
for $0 < t \le 1$\,, which implies that
$$\int_0^{\infty}t^{s-1}dt\int_M(tr \
\Ee_>(t;x,x) - \rho(x){\cdot}tr \ \Ee_{cyl}(t;x,x) - (1 -
\rho(x)){\cdot}tr \ {\tilde \Ee}(t;x,x))dx$$ is a well-defined
holomorphic function of $s$ on the whole complex plane. In
particular, we have obtained the following result

\bigskip

\begin{lem}\label{l:z1}
\begin{equation}\label{e:z3}
\z_{ \Dd_{\Pi_>}^2}(0) = \lim_{s \to
0}{\frac{1}{\G(s)}}\int_0^{\infty}t^{s-1}dt\int_M\rho(x){\cdot} tr
\ \Ee_{cyl}(t;x,x)dx \, .
\end{equation}
\end{lem}

\bigskip

\begin{proof}
Equation (\ref{e:z2}) implies the following equality
\begin{multline*}
\z_{\Dd_{\Pi_>}^2}(0)\\
= \lim_{s \to 0}{\frac{1}{\G(s)}}\int_0^{\infty}t^{s-1}dt\int_M
\Bigl(\rho(x){\cdot} tr \ \Ee_{cyl}(t;x,x) + (1 -
\rho(x)){\cdot}tr \ {\tilde \Ee}(t;x,x)\Bigr)dx \,.
\end{multline*}

It is well-known that in the case of the Dirac Laplacian on a
closed, odd-dimensional manifold, the {\it ``local"} $\z$-function
disappears (see for instance \cite{Gi95}), hence
\begin{equation}\label{e:z4}
\lim_{s \to 0}{\frac{1}{\G(s)}}\int_0^{\infty}t^{s-1}{\cdot}(1 -
\rho(x)){\cdot}tr \ {\tilde \Ee}(t;x,x)dt = 0 \, ,
\end{equation}
which gives the result.
\end{proof}

\vskip 5mm

Now let us recall that $B$ has a symmetric spectrum and its
spectral decomposition has the form
$$\{\la_n\,, \phi_n \
; \ -\la_n\,, G\phi_n\}_{n=1}^{\infty} \,.$$ The explicit
representation of the kernel $\Ee_{cyl}(t;x,y)$ with respect to
this decomposition is as follows
\begin{multline}\label{e:z5}
\Ee_{cyl}(t;(u,x),(v,y)) =
\sum_{n=1}^{\infty}{\frac{e^{-\la_n^2t}}{\sqrt{4{\pi}t}}}\{e^{-{{(u-v)
^2} \over{4t}}} - e^{-{{(u+v)^2}\over{4t}}}\}
\phi_n(x){\tensor}{\overline
{\phi_n(y)}} \\
     +
\sum_{n=1}^{\infty}{\frac{e^{-\la_n^2t}}{\sqrt{4{\pi}t}}}
\{e^{-{{(u-v)^2}\over{4t}}} + e^{-{{(u+v)^2}\over{4t}}}\}
G(x)\phi_n(x){\tensor}{\overline
{G(y)\phi_n(y)}} \\
     -
\sum_{n=1}^{\infty}\la_ne^{\la_n(u+v)}{\rm
erfc}\Bigl({{(u+v)}\over{2\sqrt{t }}} +
\la_n\sqrt{t}\Bigr)G(x)\phi_n(x){\tensor}{\overline
{G(y)\phi_n(y)}} \,,
\end{multline}
where ${\rm erfc}(r)$ is defined as in \eqref{e:erfc}:
$${\rm erfc}(r) = {\frac
2{\sqrt{\pi}}}\int_r^{\infty}e^{-\xi^2}d{\xi}\,.$$

\noi We now have an explicit representation of the integral in
(\ref{e:z3})

\begin{align}\label{e:fo1}
&\int_0^{\infty}t^{s-1}dt\int_M\rho
(x){\cdot}tr \ \Ee_{cyl}(t;x,x)dx\\
     =&
\int_0^{\infty}\rho(u)du
\int_0^{\infty}t^{s-1}2{\cdot}\Bigl\{\sum_{n=1}^{\infty}{\frac{e^{-\la_n^2t}}
{\sqrt{4{\pi}t}}}\Bigr\}dt \notag \\
    &
-
\int_0^{\infty}\int_0^{\infty}t^{s-1}\rho(u)\Bigl\{\sum_{n=1}^{\infty}\la_n
e^{2\la_nu}{\rm erfc}\Bigl({\frac{u}{\sqrt{t}}} +
\la_n\sqrt{t}\Bigr)\Bigr\}dudt \notag \\ =&
\frac{1}{2\sqrt{\pi}}{\int_0^{\infty}\rho(u)du}
\int_0^{\infty}t^{s - {3\over{2}}}Tr_Y e^{-tB^2}dt \notag
\\
&-
\int_0^{\infty}t^{s-1}dt\int_0^{\infty}\rho(u)\Bigl\{\sum_{n=1}^{\infty}
\la_ne^{2\la_nu}{\rm erfc}\Bigl({\frac{u}{\sqrt{t}}} +
\la_n\sqrt{t}\Bigr)\Bigr\}du \,.\notag
\end{align}

\noi We start with the second integral on the right side of
(\ref{e:fo1})

\begin{align*}
&\int_0^{\infty}t^{s-1}dt\int_0^{\infty}\rho(u
)\{\sum_{n=1}^{\infty}\la_ne^ {2\la_nu}{\rm
erfc}({\frac{u}{\sqrt{t}}} +
\la_n\sqrt{t})\}du \\
=& {\frac{1}{2}}\int_0^{\infty}t^{s-1}dt
\int_0^{\infty}\rho(u)\{\sum_{n=1}^{\infty}({\frac{d}{du}}e^{2u\la_n})
{\rm erfc} ({\frac{u}{\sqrt{t}}} + \la_n\sqrt{t})\}du.
\end{align*}

\noi Integration by parts leads to
\begin{align*}
& {\frac{1}{2}}\int_0^{\infty}t^{s-1}
\{[\sum_{n=1}^{\infty}\rho(u)e^{2u\la_n} {\rm
erfc}({\frac{u}{\sqrt{t}}} +
\la_n\sqrt{t})]|_0^{\infty}\}dt \\
& - {\frac{1}{2}}\int_0^{\infty}t^{s-1}dt\int_0^{\infty}\rho'(u)
\{\sum_{n=1}^{\infty}e^{2\la_nu}{\rm erfc}({\frac{u}{\sqrt{t}}}
+ \la_n\sqrt{t})\}du \\
&- {\frac{1}{2}}\int_0^{\infty}t^{s-1}dt\int_0^{\infty}\rho(u)
\{\sum_{n=1}^{\infty}e^{2\la_nu}{\rm erfc}'({\frac{u}{\sqrt{t}}} +
\la_n\sqrt{t}){\frac{1}{\sqrt{t}}}\}du \\
=& - {\frac{1}{2}}\int_0^{\infty}t^{s-1}\{\sum_{n=1}^{\infty} {\rm
erfc}(\la_n\sqrt{t})\}dt\\
& - {\frac{1}{2}}\int_0^{\infty}t^{s-1}dt\int_0^{\infty}\rho'(u)
\{\sum_{n=1}^{\infty}e^{2\la_n u}{\rm erfc}({\frac{u}{\sqrt{t}}}
+\la_n\sqrt{t})\}du \\
&+ {\frac{1}{\sqrt{\pi}}}\int_0^{\infty}t^{s-1}dt
\int_0^{\infty}\rho(u)\{\sum_ {n=1}^{\infty}e^{2\la_n
u}e^{-({{u^2}\over{t}} + 2u\la_n +
\la_n^2t)}\}{\frac{du}{\sqrt{t}}}
\end{align*}
\begin{align*}
=&
{\frac{1}{2}}\int_0^{\infty}t^{s-1}\{\sum_{n=1}^{\infty}e^{-\la_n^2t}\}dt
{\frac{2}{\sqrt{\pi}}}\int_0^{\infty}\rho(u)e^{-{{u^2}\over{t}}}{\frac{du}
{\sqrt{t}}} \\
&- {\frac{1}{2}}\int_0^{\infty}t^{s-1}\{\sum_{n=1}^{\infty} {\rm
erfc}(\la_n\sqrt{t})\}dt \\
& - {\frac{1}{2}}\int_0^{\infty}t^{s-1}dt\int_0^{\infty}\rho'(u)
\{\sum_{n=1}^{\infty}e^{2\la_nu}{\rm erfc}({\frac{u}{\sqrt{t}}} +
\la_n\sqrt{t})\}du \,.
\end{align*}
%%
%% \vfill\eject
%%
%% \noi
Finally, we have
\begin{align}\label{e:z6}
&\int_0^{\infty}t^{s-1}dt\int_0^{\infty}\rho( u)
\{\sum_{n=1}^{\infty}\la_ne^{2\la_nu}{\rm
erfc}({\frac{u}{\sqrt{t}}} +
\la_n\sqrt{t})\}du \\
=& \int_0^{\infty}t^{s-{3\over{2}}} Tr \
e^{-tB^2}dt{\cdot}{\frac{1}{2\sqrt{\pi}}}\int_0^{\infty}\rho(u)e^{-{{u^2}\over
{t}}}du \notag \\
   &  -
{\frac{1}{2}}\int_0^{\infty}t^{s-1}\{\sum_{n=1}^{\infty} {\rm
erfc}(\la_n\sqrt{t})\}dt \notag \\
   &  -
{\frac{1}{2}}\int_0^{\infty}t^{s-1}dt\int_0^{\infty}\rho'(u)
\{\sum_{n=1}^{\infty}e^{2\la_nu}{\rm erfc}({\frac{u}{\sqrt{t}}} +
\la_n\sqrt{t})\}du \,.\notag
\end{align}

\noi First, we analyze the middle term on the right side. The
following calculations hold for a single eigenvalue:
\begin{align*}
\int_0^{\infty}t^{s-1}{\rm erfc}(\la_n\sqrt{t})dt &=
{\frac{1}{s}}\int_0^{\infty}{d\over{dt}}(t^s){\rm
erfc}(\la_n\sqrt{t})dt\\
& = {\frac{t^s}{s}}{\rm erfc}(\la_n\sqrt{t})|_0^{\infty}
-{\frac{1}{s}}\int_0^{\infty}t^s{\rm erfc}'
(\la_n\sqrt{t}){\frac{\la_n}{2\sqrt{t}}}dt\\
& = -{\frac{\la_n}{2s}}\int_0^{\infty}t^{s -
{1\over{2}}}(-{\frac{2}{\sqrt{\pi}}}e^{-\la_n^2t})dt\\
& = {\frac{\la_n}{s\sqrt{\pi}}} \int_0^{\infty}t^{s -
{1\over{2}}}e^{-\la_n^2t}dt =
{\frac{\G(s+{1\over{2}})}{s\sqrt{\pi}}}\la_n^{-2s} \,.
\end{align*}

\noi It follows that for $Re(s)$ large, the middle term on the
right side of (\ref{e:z6}) is equal to
$$-
{\frac{1}{2}}\int_0^{\infty}t^{s-1} \{\sum_{n=1}^{\infty}{\rm
erfc}(\la_n\sqrt{t})\}dt = -
{\frac{1}{2}}{\frac{\G(s+{1\over{2}})}{2s\sqrt{\pi}}} \z_{B^2}(s)
\,.$$ This has a nice meromorphic extension, with simple poles, to
the whole complex plane. We rewrite (\ref{e:z6}) as
\begin{align*}
&\int_0^{\infty}t^{s-1}dt\int_0^{\infty}\rho(u)
\{\sum_{n=1}^{\infty}\la_ne^{2\la_nu}{\rm
erfc}({\frac{u}{\sqrt{t}}} + \la_n\sqrt{t})\}du \\ =&
\int_0^{\infty}t^{s-{3\over{2}}} Tr_Y
e^{-tB^2}dt{\cdot}{\frac{1}{2\sqrt{\pi}}}
\int_0^{\infty}\rho(u)e^{-{{u^2}\over{t}}}du
\\ & - {\frac{\G(s+{1\over{2}})}{4s\sqrt{\pi}}}\z_{B^2}(s) -
{\frac{1}{2}}\int_0^{\infty}t^{s-1}dt\int_0^{\infty}\rho'(u)
\{\sum_{n=1}^{\infty}e^{2\la_nu}{\rm erfc}({\frac{u}{\sqrt{t}}} +
\la_n\sqrt{t})\}du \ ,
\end{align*}
and we substitute this into (\ref{e:fo1}).

We put the final result of the computation as an independent
statement.

\bigskip

\begin{prop}\label{p:fo2}
The following equality describes the cylinder contribution to the
$\z$-function of the operator $\Dd_{\Pi_>}^2$

\begin{align}\label{e:fo2}
&\int_0^{\infty}t^{s-1}dt\int_M\rho(x){\cdot} tr \
\Ee_{cyl}(t;x,x)dx \\ =&\frac{1}{2\sqrt{\pi}}
{\int_0^{\infty}\rho(u)du}\int_0^{\infty}t^{s-{3\over{2}}} Tr_Y
e^{-tB^2}dt \notag\\ &-\{\int_0^{\infty}t^{s-{3\over{2}}} Tr_Y
e^{-tB^2}dt{\cdot}{\frac{1}{2\sqrt{\pi}}}\int_0^{\infty}\rho(u)
e^{-{{u^2}\over{t}}}du \ - \
{\frac{\G(s+{1\over{2}})}{4s\sqrt{\pi}}}\z_{B^2}(s) \notag \\ &\ \
\ \ -{\frac{1}{2}}\int_0^{\infty}t^{s-1}dt
\int_0^{\infty}\rho'(u)\{\sum_{n=1}^{\infty}e^{2\la_nu} {\rm
erfc}({\frac{u}{\sqrt{t}}} + \la_n\sqrt{t})\}du\} \,.\notag
\end{align}

\end{prop}

\bigskip

The formula (\ref{e:fo2}) is used in the study of the {\it
adiabatic decomposition} of the $\z$-determinant presented in
Section 2. We have to analyze (\ref{e:fo2}) further in order to
get information about the value of the $\z$-function at $s=0$\,.

\bigskip

\begin{lem}\label{l:hol1}
The function
\begin{equation}\label{e:hol1}
\Ff_1(s) = \int_0^{\infty}t^{s-1}dt \int_0^{\infty}\rho'(u)
\{\sum_{n=1}^{\infty}e^{2\la_nu}{\rm erfc}({\frac{u}{\sqrt{t}}} +
\la_n\sqrt{t})\}du \,
\end{equation}
is a holomorphic function on the whole complex plane.
\end{lem}

\vskip 5mm

\begin{proof}
We use the fact that $supp \ \rho' \< [{1\over{3}}\,, {2\over{3}}]
\times Y$\,, which guarantees a nice behavior of the integral with
respect to the $u$-variable since the sum over the eigenvalues is
absolutely convergent. We just have to show that $|\Ff_1(s)|$
behaves nicely with respect to $s$\,. We use the fact that ${\rm
erfc}(r) \le e^{-r^2}$ and estimate
\begin{align*}
|\Ff_1(s)| &\le \int_0^{\infty}t^{s-1}dt\int_0^{\infty}|\rho'(u)|
\{\sum_{n=1}^{\infty }e^{2\la_nu}{\rm erfc}({\frac{u}{\sqrt{t}}} +
\la_n\sqrt{t})\}du\\
& \le
%%\vfill\eject
%%
\int_0^{\infty}t^{s-1}dt\int_0^{\infty}|\rho'(u)|\{\sum_{n=1}^{\infty}
e^{-{{u^2}\over{t}} - t\la_n^2}\}du\\ & =
{\frac{1}{2}}{\cdot}\int_0^{\infty}t^{s-1}Tr_Y e^{-tB^2}dt\int_{{1\over{3}}}^{{2\over{3}}}|\rho'(u)|e^{-{{u^2}\over{t}}}du\\
&\le c_1\int_0^{\infty}t^{s-1}e^{-{{c_2}\over{t}}}Tr_Y e^{-tB^2}dt
\end{align*}
for some positive constants $c_1$\,, $c_2$ and now the Lemma
follows from the well-known asymptotics of $Tr_Y e^{-tB^2}$ as $t
\to 0$ and $t \to \infty$\,.
\end{proof}

\vskip 5mm

Now, we consider the term
\begin{equation}\label{e:hol2}
\Ff_2(s) = \int_0^{\infty}t^{s-{3\over{2}}} Tr_Y
e^{-tB^2}dt{\cdot}{\frac{1}{2\sqrt{\pi}}}
\int_0^{\infty}\rho(u)e^{-{{u^2}\over{t}}}du\, .
\end{equation}
The function $\rho(u)$ is equal to $1$ for $0 \le u \le
{1\over{3}}$ and we split the integral accordingly
\begin{align*}
\Ff_2(s) =& \int_0^{\infty}t^{s-{3\over{2}}}Tr_Y
e^{-tB^2}dt{\cdot}
{\frac{1}{2\sqrt{\pi}}}\int_0^{{1\over{3}}}e^{-{{u^2}\over{t}}}du \\
& + \int_0^{\infty}t^{s-{3\over{2}}} Tr_Y
e^{-tB^2}dt{\cdot}{\frac{1}{2\sqrt{\pi}}}
\int_{{1\over{3}}}^{{2\over{3}}}\rho(u)e^{-{{u^2}\over{t}}}du \,.
\end{align*}
Let us observe that
\begin{align*}
\int_0^{{1\over{3}}}e^{-{{u^2}\over{t}}}du &=
\sqrt{t}{\cdot}\int_0^{1\over{3\sqrt{t}}}e^{-y^2}dy =
\sqrt{t}{\cdot}\int_0^{\infty}e^{-y^2}dy -
\sqrt{t}{\cdot}\int_{1\over{3\sqrt{t}}}^{\infty}e^{-y^2}dy\\
& = {\frac{\sqrt{\pi}}{2}}\sqrt{t} -
{\frac{\sqrt{\pi}}{2}}\sqrt{t}{\cdot} {\rm
erfc}({1\over{3\sqrt{t}}}) \ \,,
\end{align*}
which allows us to represent $\Ff_2(s)$ in the following form
\begin{align}\label{e:ho13}
\Ff_2(s) =&{\frac{1}{4}}\int_0^{\infty}t^{s-1} Tr_Y e^{-tB^2}dt -
{\frac{1}{4}}\int_0^{\infty}t^{s-1}
Tr_Y e^{-tB^2}{\cdot}{\rm erfc}({1\over{3\sqrt{t}}})dt \\
& + \int_0^{\infty}t^{s-{3\over{2}}}Tr_Y e^{-tB^2}dt
{\cdot}{\frac{1}{2\sqrt{\pi}}} \int_{{1\over{3}}}^{{2\over{3}}}
\rho(u)e^{-{{u^2}\over{t}}}du \ \,.\notag
\end{align}

\noi The middle term on the right side of the above equality is
again holomorphic on the whole complex plane due to the inequality
$$\Bigl|\int_0^{\infty}t^{s-1}Tr_Y e^{-tB^2}{\cdot}{\rm
erfc} \bigl({\frac{1}{3\sqrt{t}}}\bigr)dt\Bigr| \le c
\int_0^{\infty}t^{s-1}Tr_Y e^{-tB^2}{\cdot}e^{-{1\over{9t}}}dt \ \
.$$ We estimate the last term on the right side of (\ref{e:ho13})
in the same way to show that it is a holomorphic function of $s$
as well. Finally, we evaluate the $\z$-function at $s = 0$, using
Lemma \ref{l:z1}:
\begin{align*}
     \z_{\Dd_{\Pi_>}^2}(0) &= \lim_{s \to
0}{\frac{1}{\G(s)}} \int_0^{\infty}t^{s-1}dt\int_M\rho(x){\cdot}tr
\
\Ee_{cyl}(t;x,x)dx \\
&=\lim_{s \to 0}{\frac{1}{\G(s)}}
\{\frac{1}{2\sqrt{\pi}}\int_0^{\infty}{\rho(u)du}\int_0^
{\infty}t^{s -
{3\over{2}}}Tr_Y e^{-tB^2}dt\\
&\qquad\qquad \quad - \Ff_2(s)  +
{\frac{\G(s+{1\over{2}})}{4s\sqrt{\pi}}}\z_{B^2}(s) +
{\frac{1}{2}}\Ff_1(s)\}\\
& = \lim_{s \to 0}s{\cdot}
\{\frac{1}{2\sqrt{\pi}}\int_0^{\infty}{\rho(u)du}\int_0^{\infty}t
^{s -
{3\over{2}}}Tr_Y e^{-tB^2}dt\\
&\qquad\qquad - \Ff_2(s) +
{\frac{\G(s+{1\over{2}})}{4s\sqrt{\pi}}}\z_{B^2}(s) +
{\frac{1}{2}}\Ff_1(s)\} \\
& = \lim_{s \to 0}s{\cdot}
\{\frac{1}{2\sqrt{\pi}}\int_0^{\infty}{\rho(u)du}\int_0^{\infty}t
^{s -
{3\over{2}}}Tr_Y e^{-tB^2}dt\\
&\qquad\qquad - {\frac{1}{4}}\int_0^{\infty}t^{s-1}Tr_Y
e^{-tB^2}dt +
{\frac{\G(s+{{1}\over{2}})}{4s\sqrt{\pi}}}\z_{B^2}(s)\}\\
& = \lim_{s \to
0}s{\cdot}\frac{1}{2{\sqrt{\pi}}}\int_0^{\infty}{\rho(u)du}
\int_0^{\infty}t
^{s - {3\over{2}}}Tr_Y e^{-tB^2}dt \\
&\ \ + \lim_{s \to
0}\{s{\frac{\G(s+{1\over{2}})}{4s\sqrt{\pi}}}\z_{B^2}(s) -
s{\frac{1}{4}}\int_0^{\infty}t^{s-1}Tr_Y e^{-tB^2}dt\}\\
& = 0 + \{{\frac{1}{4}}\z_{B^2}(0) - {\frac{1}{4}}\z_{B^2}(0)\} =
0 \,.
\end{align*}

\bigskip

The situation is not different in the case of non-invertible
$\Dd_{\Pi_>}$ . We have

$$\z_{\Dd_{\Pi_>}^2}(0) = \lim_{s
\to 0}{\frac{1}{\G(s)}}\int_0^{\infty}t^{s-1} (Tr \
e^{-t\Dd_{\Pi_>}^2} - dim \ ker \Dd_{\Pi_>}^2)dt \ \ ,$$

\bigskip
\noi where the dimension of the kernel is present in order to make
the integral $\int_1^{\infty}$ convergent. Now we have

\begin{align*}
\z_{\Dd_{\Pi_>}^2}(0) =& \lim_{s\to 0}{\frac{1}{\G(s)}}
\int_0^1t^{s-1}(Tr \ e^{-t\Dd_{\Pi_>}^2} -
dim \ ker \ \Dd_{\Pi_>}^2)dt \\
=&\lim_{s \to 0}({\frac{1}{\G(s)}}\int_0^1t^{s-1}Tr \
e^{-t\Dd_{\Pi_>}^2}dt) -
dim \ ker \ \Dd_{\Pi_>}^2 \\
=&\lim_{s \to 0}({\frac{1}{\G(s)}}
\int_0^{\infty}t^{s-1}dt\int_M\rho(x){\cdot}tr \
\Ee_{cyl}(t;x,x)dx) -
dim \ ker \ \Dd_{\Pi_>}^2 \\
=& - dim \ ker \ \Dd_{\Pi_>}^2 \ \ .
\end{align*}

\bigskip
We also do not have problem with the case $ker \ B \ne \{0\}$ .
The {\it  Cobordism Theorem} for the Dirac operators (see for
instance \cite{BoWo93}) implies the existence of the involution
$$\s : ker \ B \to ker \ B \ \ ,$$
such that $G\s = - {\s}G$ . Let $\pi_{\s} : ker \ B \to ker \ B$
denote orthogonal projection onto $+1$-eigenspace of $\s$ . The
orthogonal projection

\noi $\Pi_{\s} = \Pi_> + \pi_{\s}$ is an element of
$Gr_{\infty}^*(\Dd)$ and we can repeat our computations to obtain
$$\z_{\Dd_{\Pi_{\s}}^2}(0) =
- dim \ ker \ \Dd_{\Pi_{\s}}^2 \ \ .$$

\bigskip

Finally the result for arbitrary element $P \in
Gr_{\infty}^*(\Dd)$ follows from the existence of a positive
constant $c > 0$ , such that for any $0 < t < 1$

$$|Tr \ e^{-t\Dd_P^2} - Tr \
e^{-t\Dd_{\Pi_{\s}}^2}| < c\sqrt{t} \ \ .$$

\bigskip

\noi This result is stated as Theorem 3.2 in \cite{KPW99}. The
proof consists of a straightforward computation and the details
are presented in \cite{KPW99}. The idea is easy to understand. It
was explained in Section 1 of \cite{KPW99}, that $\Dd_P^2$ is
unitarily equivalent to the operator of the form $\Dd_{\Pi_{\s}}^2
+ \Kk$ , where $\Kk : L^2(M;S) \to L^2(M;S)$ is a bounded
operator, with kernel $\Kk(x,y)$ supported in $N = [0,1] \times Y$
. Moreover, $\Kk(x,y)$ is smoothing in $Y$-direction. By the
Duhamel's Principle we have

$$Tr \ e^{-t\Dd_P^2} - Tr \ e^{-t\Dd_{\Pi_{\s}}^2} = -
Tr \int_0^te^{-s\Dd_P^2}\Kk e^{-(t-s)\Dd_{\Pi_{\s}}^2}ds \ \ .$$

\bigskip
\noi The expression on the right side can be written as the
series, where each next term has the better behavior with respect
to $t$ , than the previous one. The first term is

$$- Tr
\int_0^te^{-s\Dd_{\Pi_{\s}}^2}\Kk e^{-(t-s)\Dd_{\Pi_{\s}}^2}ds = -
\int_0^tTr \ \Kk e^{-t\Dd_{\Pi_{\s}}^2} = -t{\cdot}Tr \ \Kk
e^{-t\Dd_{\Pi_{\s}}^2} \ \ .$$

\bigskip
\noi Now the kernel of the operator $\Kk$ is smoothing in the
$Y$-direction, hence the only singularity left is in the normal
direction and we obtain

$$|Tr \
e^{-t\Dd_P^2} - Tr \ e^{-t\Dd_{\Pi_{\s}}^2}| \sim_{t\to 0} t|Tr \
\Kk e^{-t\Dd_{\Pi_{\s}}^2}| \le t{\cdot}{c\over{\sqrt{t}}} \le
c{\cdot}\sqrt{t}$$

\bigskip

\noi (we refer to \cite{KPW99} for the detailed presentation). It
follows that

$$|\lim_{s \to
0}{\frac{1}{\G(s)}}\int_0^{\infty}t^{s-1}(Tr \ e^{-t\Dd_P^2} - Tr
\ e^{-t\Dd_{\Pi_{\s}}^2})dt| \le$$

$$\lim_{s \to 0}{\frac{1}{\G(s)}}\int_0^1t^{s-1}
|Tr \ e^{-t\Dd_P^2} - Tr \ e^{-t\Dd_{\Pi_{\s}}^2}|dt \le$$

$$c{\cdot}\lim_{s \to 0}s{\cdot}\int_0^1t^{s-{1\over{2}}}ds = 0 \ \ ,$$

\bigskip

and as a result we have

$$\z_{\Dd_P^2}(0) - \z_{\Dd_{\Pi_{\s}}^2}(0) =
dim \ ker \ \Dd_{\Pi_{\s}} - dim \ ker \Dd_P \ \ .$$
\bigskip

This ends the proof of the Proposition \ref{p:zeta}.

\end{appendix}

\end{document}